# On the Blind Source Separation of Human Electroencephalogram by Approximate Joint Diagonalization of Second Order Statistics

Marco Congedo[*], Cédric Gouy-Pailler, Christian Jutten

GIPSA-lab (Grenoble Image Parole Signaux Automatique),
UMR5216 :  CNRS (Centre National de la Recherche Scientifique) - Université Joseph Fourier - Université Pierre Mendès-France - Université Stendhal - INPG (Institut Polytechnique de Grenoble)
* Corresponding author: CNRS, GIPSA-lab. 46, avenue Félix Viallet 38031 GRENOBLE Cedex
Tel: + 33 (0)4 7657 4352; Fax: + 33 (0)4 7657 4790
e-mail: Marco.Congedo@Gmail.com

Review Paper, Clinical Neurophysiology (2008), 119,  2677–2686

**Acknowledgements:**
This Research has been partially supported by the French National Research Agency (ANR) within the National Network for Software Technologies (RNTL), project Open-ViBE ("Open Platform for Virtual Brain Environments"), grant # ANR05RNTL01601, by the European COST Action B27 "Electric Neuronal Oscillations and Cognition". During the period of the research the first author has been partially supported by Nova Tech EEG, Inc., Knoxville, TN and the second author by the French Ministry of Defense (DGA).The authors wish to express their gratitude to Berrie Gerrits, Bertrand Rivet, Reza Sameni, Leslie Sherlin, Antoine Souloumiac and the anonymous reviewers for their helpful comments and suggestions about this paper.

**Abstract**

Over the last ten years blind source separation (BSS) has become a prominent processing tool in the study of human electroencephalography (EEG). Without relying on head modeling BSS aims at estimating both the waveform and the scalp spatial pattern of the intracranial dipolar current responsible of the observed EEG. In this review we begin by placing the BSS linear instantaneous model of EEG within the framework of brain volume conduction theory. We then review the concept and current practice of BSS based on second-order statistics (SOS) and on higher-order statistics (HOS), the latter better known as independent component analysis (ICA). Using neurophysiological knowledge we consider the fitness of SOS-based and HOS-based methods for the extraction of spontaneous and induced EEG and their separation from extra-cranial artifacts. We then illustrate a general BSS scheme operating in the time-frequency domain using SOS only. The scheme readily extends to further data expansions in order to capture experimental source of variations as well. A simple and efficient implementation based on the approximate joint diagonalization of Fourier cospectral matrices is described (AJDC). We conclude discussing useful aspects of BSS analysis of EEG, including its assumptions and limitations.

**Introduction**

Recent studies on human electroencephalogram (EEG) are based on the theory of brain volume conduction. It is well established that the generators of brain electric fields recordable from the scalp are macroscopic post-synaptic potentials created by assemblies of pyramidal cells of the neocortex (Speckmann and Elger, 2005). Pyramidal cells are aligned and oriented perpendicularly to the cortical surface. Their synchrony is possible thanks to a dense net of local horizontal connections (mostly <1mm). At recording distances larger than about three/four times the diameter of the synchronized assemblies the resulting potential behaves as if it were produced by electric dipoles; all higher terms of the multipole expansion vanish and we obtain the often invoked dipole approximation (Lopes Da Silva and Van Rotterdam, 2005; Nunez and Srinivasan, 2006, Ch. 3). Three physical phenomena are important for the arguments we advocate in this study. First, unless dipoles are moving there is no appreciable delay in the scalp sensor measurement (Lopes da Silva and Van Rotterdam, 2005). Second, in brain electric fields there is no appreciable electro-magnetic coupling (magnetic induction) in the frequencies up to about 1MHz, thus the quasi-static approximation of Maxwell equations holds throughout the spectrum of interest (Nunez and Srinivasan, 2006, p. 535-540). Finally, for source oscillations below 40Hz it has been verified experimentally that capacitive effects are also negligible, implying that potential difference is in phase with the corresponding generator (Nunez and Srinivasan, 2006, p. 61). These phenomena strongly support the *superposition principle*, according to which the relation between neocortical dipolar fields and scalp potentials may be approximated by a system of linear equations (Sarvas, 1987). Whether this is a great simplification, we need to keep in mind that it does not hold true for all cerebral phenomena. Rather, it does at the macroscopic spatial scale we are interested in here.

A common approach to the study of human EEG is to describe patterns in space and time and link empirical findings with anatomical and physiological knowledge. The problem is characterized by high temporal resolution (about 1ms) and low spatial resolution (several $cm^3$). For example, it has been estimated that without time averaging about 60 million contiguous neurons must be synchronously active as to produce observable scalp potentials (Nunez and Srinivasan, 2006, p. 21).

Such a cluster would realistically extend over several cm$^2$ of cortical gyral surface, whereas disentangling fields emitted by cortical functional units may require much higher precision. Because of volume conduction, scalp EEG potentials describe a *mixture* of the fields emitted by several dipoles extending over large cortical areas. Practically, in order to improve the spatial resolution it is often necessary to trade in the temporal one operating some form of temporal averaging. In summary, the path followed by much of current EEG research is to "isolate" in space and time the generators of the observed EEG as much as possible, counteracting the mixing caused by volume conduction and maximizing the signal-to-noise ratio (SNR).

Over the years we have assisted to the development of several classes of methods to improve the spatial specificity. Those include, among others, surface and cortical Laplacian (Nunez and Srinivasan, 2006), equivalent dipole fitting (Mosher et al., 1992) and distributed minimum norm (model-driven) or minimum variance (data-driven) inverse solutions (Greenblatt et al., 2005; Lopes da Silva, 2004). Targeted attempts include sparsification approaches (Gorodnitsky et al., 1995; Cotter et al., 2005) and spatial filters known as beamformers (Rodríguez-Rivera et al., 2006; Congedo, 2006). Surface Laplacian methods apply a spatial high-pass filtering to the scalp potential by estimating their second spatial derivative. They tend to overemphasize high spatial frequency and radial (to the scalp surface) dipolar fields. Inverse solutions seek *source localization* in a chosen solution space and rely on geometrical models of the head tissue. Unfortunately, the accurate description of EEG volume conduction is complicated by inhomogeneity (resistivity varies with type of tissue) and anisotropy (resistivity varies in different directions); therefore source localization methods are inevitably undermined by geometrical modeling error.

Another approach that persists in EEG literature is blind source separation (BSS). First studied in our laboratory during the first half of the 80's (Ans et al., 1985; Hérault and Jutten, 1986) BSS has enjoyed considerable interest worldwide only a decade later, inspired by the seminal papers of Jutten and Hérault (1991), Comon (1994) and Bell and Sejnowski (1995). BSS has today greatly expanded encompassing a wide range of engineering applications such as speech enhancement, image processing, geophysical data analysis, wireless communication and biological signal analysis (Hyvärinen et al., 2001; Cichocki and Amari, 2002; Choi et al., 2005). Such ubiquity springs from the "blind" nature of the BSS problem formulation: no knowledge of volume conduction or of source waveform is assumed. The problem may be attacked from several perspectives; several hundred BSS algorithms have been proposed over the last 20 years with more added on every year. Typically, such methods are based on the cancellation of second order statistics (SOS) and/or of higher (than two) order statistics (HOS). Their commonality resides in the assumption of a certain degree of source *spatial independence*, which is precisely modeled by the cancellation of those statistics. Both HOS and SOS have been employed with success in EEG. They are today established for denoising/artifact rejection (Vigário, 1997; Jung et al., 2000; Vorobyov and Cichocki, 2002; Iriarte et al., 2003; Joyce et al., 2004; Kierkels at al., 2006; Fitzgibbon et al., 2007; Frank and Frishkoff, 2007; Halder at al, 2007; Phlypo et al., 2007; Romero et al., 2008; Crespo-Garcia et al., 2008), improving brain computer interfaces (Qin et al., 2005; Serby et al., 2005; Wang and James, 2007; Dat and Guan, 2007; Kachenoura et al., 2008) and for increasing the SNR of single-trial time-locked responses (Cao et al., 2002; Sander et al., 2005; Lemm et al., 2006; Tang et al., 2006; Guimaraes et al., 2007; Zeman et al., 2007). Yet, it appears that only four of the many existing algorithms have repeatedly occurred in EEG literature. They are known as FastICA (Hyvärinen, 1999), JADE (Cardoso and Souloumiac, 1993), InfoMax (Bell and Bejnowsky, 1995) and SOBI (Belouchrani et al., 1997). FastICA, InfoMax and JADE are ICA (HOS) methods, while SOBI is a SOS method. JADE and SOBI are solved by *approximate joint diagonalization* (Cardoso and Souloumiac, 1993; Pham, 2001 b; Yeredor, 2002; Ziehe et al., 2004; Vollgraf and Obermayer, 2006; Li and Zhang, 2007; Fadaili et al, 2007; Dégerine and Kane, 2007), a powerful algebraic tool which allows promising extensions that we will consider in this study.

**The BSS problem for the brain**

For $N$ scalp sensors and $M \leq N$ EEG dipolar fields with fixed location and orientation in the analyzed time interval, the linear BSS model simply states the superposition principle discussed above, i.e.,

$$v(t) = As(t) + \eta(t), \qquad (1.0)$$

where $v(t) \in \mathbb{R}^N$ is the *sensor measurement vector*, $A \in \mathbb{R}^{N \cdot M}$ is a time-invariant full column rank *mixing matrix*, $s(t) \in \mathbb{R}^M$ holds the time-course of the source components and $\eta(t) \in \mathbb{R}^N$ is additive noise, temporally white, possibly uncorrelated to $s(t)$ and with spatially uncorrelated components. Our source estimation is given by

$$\hat{s}(t) = \hat{B}v(t), \qquad (1.1)$$

where $B \in \mathbb{R}^{M \cdot N}$ is called the *demixing* or *separating matrix*. Hereafter the caret indicates a statistical estimation. Although this is the classical BSS model we need a few clarifications for the EEG case: first, by $\eta(t)$ we model *instrumental* noise only. In the following we drop the $\eta(t)$ term because the instrumental (and quantization) noise of modern EEG equipment is typically low (<1μV). On the other hand, *biological* noise (extra-cerebral artifacts such as eye movements and facial muscle contractions) and *environmental* noise (external electromagnetic interference) may obey a mixing process as well, thus they are generally modeled as components of $s(t)$, along with cerebral ones. Notice that while biological and environmental noise can be identified as separated components of $s(t)$, hence removed, source estimation will be affected by the underlying cerebral *background noise* propagating with the same coefficients as the signal (Belouchrani and Amin, 1998). Second, the assumption of time-invariance of the mixing process in (1.0) must apply only locally. The demixing matrix is assumed fixed for a given temporal interval, but may be allowed to change (slowly) across successive intervals (Pham 2001 a; Li et al., 2006). Such a model allows changes in the location and orientation of dipole layers over time. The assumptions underlying model (1.0) are crucial for the success of the source separation, thus will be reconsidered in more details in the discussion.

**A suitable class of solutions to the brain BSS problem**

To tackle problem (1.1) assuming knowledge of sensor measurement only we need to reduce the number of admissible solutions. In this paper we are interested in weak restrictions converging toward condition

$$\hat{s}(t) = Gs(t), \qquad (1.2)$$

where $s(t)$ holds the time-course of the true (unknown) source processes and the *system matrix*

$$G = \hat{B}A \approx \Lambda P \qquad (1.3)$$

approximates a signed scaling (a diagonal matrix $\Lambda$) and raw permutation ($P$). Equation (1.2) is obtained substituting (1.0) in (1.1) ignoring the noise term in the former. Whether condition (1.2) may be satisfied is a problem of *identifiability*, which establish the theoretical ground of BSS theory (Tong et al., 1990; Tong et al., 1991 a, b; Tong et al., 1993; Cardoso, 1998 a; Pham and Cardoso, 2001; Pham, 2002; Theis, 2004). In turn, matching condition (1.2) implies that we can recover faithfully the source *waveform* out of a *scale* (including sign) and *permutation* (order) indeterminacy. The idea suits EEG well, since the waveform bears meaningful physiological and clinical information. Notice the correspondence between the $m^{th}$ source, its *separating vector* ($m^{th}$ row of $\hat{B}$) and its *scalp spatial pattern* (mixing vector), given by the $m^{th}$ *column* of $\hat{A} = \hat{B}^+$. Hereafter superscript + indicates the Moore-Penrose pseudo-inverse. The mono-dimensionality of those vectors and their sign/energy indeterminacy implies the explicit modeling of the orientation and localization parameters of the $m^{th}$ source, but not its moment. This is also the case of inverse solutions with good source localization performance (Greenblatt et al., 2005). Nonetheless, we can evaluate the relative energy of each source sorting them by decreasing value of *explained variance*, as illustrated in appendix (D). Moreover, linearity allows switching back from the source space into the sensor space. Substituting (1.1) into (1.0) and dropping the noise term in the latter yields *BSS filtering*

$$v'(t) = \hat{A}R\hat{s}(t) = \hat{A}R\hat{B}v(t), \qquad (1.4)$$

where $R$ is a diagonal matrix with $m^{th}$ diagonal element equal to 1 if the $m^{th}$ component is to be retained and equal to 0 if it is to be removed.

**Different approaches for solving the source separation problem**

It has been known for a long time that in general the BSS problem cannot be solved for sources that are Gaussian, independent and identically distributed (iid) (Darmois, 1953). The iid condition implies that each sample of the source components is statistically independent from the others and that they all follow the same probability distribution. Therefore, in order to solve the BSS problem the sources must be either (1) possibly iid, but non-Gaussian or (2) not iid. In case (1), one assumes that at most one source is Gaussian and that they are all mutually statistically independent. The mutual independence assumption (spatial independence of all pair-wise sources) should not be confused with the iid condition (temporal independence of successive samples within each source process). Actually, the iid condition implies that no temporal information is used, thus the method is efficient regardless the temporal dependence of sources. Those methods are known as *independent component analysis* (ICA) (Jutten and Hérault, 1991; Comon, 1994, 1999; Hyvärinen et al., 2001). ICA requires higher order statistics (HOS), explaining why it may succeed only if at most one source has Gaussian distribution: in fact Gaussian distributions are fully defined by their statistics up to the second order (SOS). The idea of (2) is to break the non-Gaussianity assumption. This can be done by assuming that source components are all pair-wise uncorrelated and that either (a) within each source component the successive samples are temporally correlated[1], (Tong et al., 1990; Molgedey and Schuster, 1994; Belouchrani et al., 1997; Ziehe and Müller, 1998) or (b) samples in successive time intervals do not have the same statistical distribution, i.e., they are non stationary (Matsuoka et al., 1995; Souloumiac, 1995; Choi and Cichocki, 2000; Pham and Cardoso, 2001; Choi et al., 2002). Provided that source components have non-proportional spectra or the time courses of their variance (energy) vary differently, one can show that second order statistics are sufficient for solving the source separation problem. Since second order statistics are sufficient, the method is able to separate also Gaussian sources, contrary to ICA.

In summary, ICA methods require a higher independence assumption (HOS independence), while SOS methods rely on a weaker uncorrelation assumption (SOS independence) coupled to the assumption that source components display unique spectral density signature (characteristic source coloration) and/or unique source energy variation signature (characteristic non-stationarity). We will see that different source energy in different experimental conditions is also a sufficient additional condition allowing separation using SOS. If these assumptions are fulfilled the separating matrix can be identified uniquely, thus source can be recovered *regardless the true mixing process* (uniform performance property: see for example Cardoso, 1998 a) and regardless the iid condition for ICA or the Gaussian condition for SOS BSS.

---

[1] Such processes are called *colored*, in opposition to iid processes, which are called *white*.

A recent trend in the literature on EEG source separation is the design of algorithms blending spatial independence with other specific assumptions. Such an approach is called *semi-blind source separation*. Different priors have been introduced in the cost function pursuing spatial independence. In a Bayesian framework Roberts (1998) introduces priors on signal distributions, but those seems little useful in EEG since typically the distributions are not known a-priori. Temporal constraints have been introduced using ad hoc reference signals based on experimental stimulation or visual inspection of the observed mixtures (James and Gibson, 2003; Lu and Rajapakse, 2005). Zhang (2008) describes a method to extract reference signals from the observed waveforms designed to work even when the interesting source signal is not visible in the observed measurement. Spatial constraints have been employed to recover sources of interest with known spatial topographies (Hesse and James, 2006; Ille, Berg and Scherg 2002). Spectral constraints have been also introduced to recover sources with known spectral content (James and Hesse, 2005; Wang and James, 2007; Barbati et al, 2008). Finally, a flexible semi-blind ICA approach where one can incorporate priors of experimental or physiological origin as well has been proposed by Barbati et al. (2006). All BSS methods cited so far rely on the spatial independence assumption. Other methods exist, relaxing definitely this assumption and imposing instead positivity of sources and mixtures (Lee and Seung, 1999, 2001) or source sparsity (Gribonval and Lesage, 2006; Li et al., 2006). These methods will not be considered in this study despite their theoretical interest since so far their use with EEG has been marginal.

**SOS vs. HOS: statistical considerations**

Joyce et al. (2004) reports that successful separation of EEG data can be achieved using as few as 100 data points using SOBI (SOS) and 1000 or more using ICA (HOS) algorithms. This is a known advantage of SOS-based BSS methods; an higher statistical efficiency allows performing BSS on shorter time intervals, which is a safe strategy to prevent serious departures from the linear instantaneous model assumptions (see discussion). It has also been suggested that SOS estimations are more robust with respect to noise as compared to HOS estimations (Belouchrani et al, 1997; Joyce et al, 2004). We notice that this is not necessarily true. For instance, the estimation of kurtosis is unaffected by white Gaussian noise, whereas this is not true for SOS estimations. However, SOS estimations are expected to be more robust with respect to outliers, in that their influence in ensemble average estimations is magnified by elevation to the second power, whereas for HOS it is magnified by elevation to the third and/or fourth power.

**SOS vs. HOS: neurophysiological considerations**

*The hypothesis of spatial independence*

Human neocortex is a prodigious net of local and global interconnections. There are about as many neurons ($10^{10}$) as cortico-cortical fibers connecting them in the 1-15 cm range (Nunez and Srinivasan, 2006, p. 7). Dense and sometimes distributed connections exist between the neocortex and sub-cortical structures as well. Therefore, one may ask if assuming independent time course of cortical cell assemblies is reasonable. It has been speculated that forcing independence of the BSS output may result in spurious source components (Li et al., 2006). A study on induced visual gamma activity has also questioned the source non-Gaussianity assumption as required by ICA (Barbati et al., 2008). In practice, the BSS output is never exactly independent, but just as independent as possible and this may explain why BSS is useful with EEG. HOS and SOS BSS may be both conceived as spatial filters minimizing the dependence of the observed mixtures. For EEG data this is an effective way to counteract the effect of volume conduction. In fact, we have seen that the brain tissue behaves approximately as a linear conductor, thus observed potentials (mixtures) are more dependent than the generating dipolar fields. The fundamental difference between HOS and SOS BSS springs from the kind of statistical information they try to extract from the data. Thus, the critical question is how well the respective statistics inform about the actual source process, question to which we now turn.

*Extra-cerebral Artifacts*

An important problem with EEG recording is the contamination of extra-cerebral artifacts. The most common artifacts are electric signals produced by the eyes and muscles of the face, jaw and neck. Artifacts have characteristic spatial, coloration and non stationary signatures, well distinct from EEG spontaneous activity (Ille et al., 2002; Lopes da Silva, 2005 a). Still, their separation and removal has remained a difficult task. Eye blinks are attributed to change of conductance due to the moving eyelid on the cornea. They generate a potential peak which amplitude can be one order of magnitude stronger than the EEG. Eye movements, both volitional movements and saccades, are attributed to the cornea-retina dipolar field. Typically, vertical and horizontal movements are monitored by means of electrodes positioned above and next to the eye, providing an orthogonal space for movements in all possible direction. On the other hand, contamination by facial, jaw and neck muscles typically manifests as persistent low-voltage high-frequency signals (>20Hz: Whitham et al., 2007) with focused spatial distribution and rapidly decaying autocorrelation function.

The use of coloration has been suggested for separating EEG from eye movements (Joyce et al., 2004). In most comparative studies (Kierkels at al., 2006; Fitzgibbon et al., 2007; Halder at al, 2007; Romero et al., 2008) coloration have been found superior than non-Gaussianity for eye movements removal, but Phlypo et al. (2007) found inconclusive results (slightly favoring HOS) and Frank and Friskhoff (2007) concluded clearly in favor of HOS. We notice that the linear instantaneous model (1.0) does *not* admit rotating dipoles. Furthermore, the movements of the two eyes are extremely correlated, forming two spatially dependent dipoles. Insomuch, in principle linear BSS methods cannot resolve them. In effect, a method assuming instantaneous mixing will try to model an "average" big dipole covering both eyes, or in the middle of them and will try to explain the rotating orientation by two or three orthogonal components. This raises the problem of "misallocation of variance" and will result in sub-optimal separation (Frank and Frishkoff, 2007). Based on these considerations we believe that a general BSS approach based on the spatial independence assumption is not optimal for separating EEG from eye movements, regardless the use of SOS or HOS. To accomplish this task semi-blind approaches appear more promising (e.g., Ille, Berg & Scherg 2002).

On the other hand, the separation of muscle contamination appears easier to treat. Generally both SOS and HOS methods perform correctly with them (Crespo-Garcia et al., 2008), yet research on more targeted approaches is in progress (Gasser et al., 2005; De Clercq et al., 2006).

*Spontaneous and induced EEG*

Observed potentials are the summation of post-synaptic potentials over large cortical areas caused by trains of action potentials carried by afferent fibers. The action potentials come in trains/rest periods, resulting in sinusoidal oscillations of the scalp potentials, with negative shifts during the train discharges and positive shifts during rest. The periodicity of trains/rest periods are deemed responsible for high-amplitude EEG rhythms (oscillations) up to about 12Hz, whereas higher frequency (>12Hz) low-amplitude rhythms may result from sustained (tonic) afferent discharges (Speckmann and Elegr, 2005). There is no doubt that an important portion of spontaneous EEG activity is rhythmic, whence strongly colored (Niedermeyer, 2005 a; Steriade, 2005; Buzsáki, 2006, Ch. 6, 7). Some rhythmic waves come in more or less short bursts. Typical examples are sleep spindles (7-14Hz) (Niedermeyer, 2005 b; Steriade, 2005), frontal Theta (4-7Hz) and Beta (13-35Hz) waves (Niedermeyer, 2005 a). Others are more sustained, as it is the case for slow Delta (1-2Hz) waves during deep sleep stages III and IV (Niedermeyer, 2005 b), the Rolandic Mu rhythms (around 10Hz and 20Hz) and posterior Alpha rhythms (8-12Hz) (Niedermeyer, 2005 a). In all cases brain electric oscillations are not ever-lasting and one can always define time intervals when rhythmic activity is present and others when it is absent or substantially reduced. Such intervals may be precisely defined based on known reactivity properties of the rhythms. For example, in event-related synchronization/desynchronization (ERD/ERS: Pfurtscheller and Lopes da Silva, 2004), which are time locked, but not phase locked increases/decreases of the oscillating energy (Steriade, 2005), intervals may be defined before and after event onset. On the other hand event-related potentials (ERP: Lopes Da Silva, 2005 b), which are both time-locked and phase-locked can be further partitioned in several successive intervals comprising the different peaks. Such source energy variation signatures can be modeled precisely by SOS, as we will specify.

*Transients*

Another class of brain electric phenomena comprises transient waves such as spikes, sharp waves and spike-wave complexes in epileptic disorder (Niedermeyer, 2005 c), vertex waves during sleep (Niedermeyer, 2005 b), etc. Transients are characterized by abrupt and sometimes large potential shifts. In general, those are not naturally characterized by coloration, unless they results from the superposition of several continuous colored waves. This is the case, for example, of K-complexes observed during sleep (Niedermeyer, 2005 b), which are a superposition of a slow wave (<1Hz) and a Delta wave (1-4Hz) (Steriade, 2005). Nonetheless, transient activities are by definition spaced by intervals of inactivity, hence the difference between the energy in their active and inactive intervals (non stationarity) may be captured adequately by SOS statistics. However, due to their possible highly non-Gaussian nature, this kind of phenomena is naturally modeled by HOS statistics.

In summary, it appears that a wide variety of spontaneous and induced EEG phenomena are captured appropriately by SOS statistics, however for transient activity HOS may be better candidates. So far SOS methods applied to EEG have concentrated mainly on coloration (e.g., SOBI). The validity of the coloration assumption for recovering actual EEG dipolar fields has received experimental support (Tang et al., 2004; Sutherland and Tang, 2006; Van Der Loo et al., 2007). We have contended that source energy variation over time is a ubiquitous property of EEG and it should be exploited besides coloration. This is the focus of SOS time-frequency approaches, which are well established in other technical fields (Belouchrani and Amin, 1998; Pham 2002; Choi et al., 2002; Bousbia-Salah et al., 2003). Here we pursue further this path in the context of experimental and clinical EEG data.

**Approximate joint diagonalization**

The class of SOS BSS methods we are considering is consistently solved by *approximate joint diagonalization* algorithms (Cardoso and Souloumiac, 1993; Pham, 2001 b; Yeredor, 2002; Ziehe et al., 2004; Vollgraf and Obermayer, 2006; Li and Zhang, 2007; Fadaili et al, 2007; Dégerine and Kane, 2007). Given a set of matrices {$Q_1, Q_2,...$ }, the AJD seeks a matrix $\hat{B}$ such that the products $\hat{B}Q_1\hat{B}^T$, $\hat{B}Q_2\hat{B}^T$, ... are as diagonal as possible (subscript "$T$" indicates matrix transposition) . Given an appropriate choice of the *diagonalization set* {$Q_1, Q_2,...$ } such matrix $\hat{B}$ is indeed an estimation of the separating matrix in (1.1) and one obtain an estimate of the mixing matrix as $\hat{A} = \hat{B}^+$. Matrices in {$Q_1, Q_2,...$ } are chosen so as to hold in the off-diagonal entries statistics describing some form of *dependence* among the sensor measurement channels; then the AJD will vanish those terms resulting in linear combination vectors (the rows of $\hat{B}$) extracting "independent" components from the observed mixture via (1.1). More particularly, the joint diagonalization is applied on matrices that *change* according to the assumptions about the source. They are those changes, when available, that provide enough information to solve the BSS problem. If the source process is assumed to be colored, one may consider lagged covariance matrices of signals. If the source process is assumed to be non stationary between blocks of data, one may consider covariance matrices estimated on different time windows. In both situations, provided that source spectra are non proportional (colored sources) or source energy varies differently (non stationary sources), the additional matrices add information (in fact, equations) sufficient for estimating all the parameters of the separating system. If the source is both colored and non stationary, one can use a set of both kinds of matrices, as we will illustrate.

The aforementioned popular JADE and SOBI algorithms are based on AJD and this is the case for many other BSS algorithms (for a review see Theis and Inouye, 2006). One advantage of AJD algorithms is that they execute fast and do not require setting parameters for convergence. In particular, the algorithms by Pham (2001 b) and by Ziehe et al. (2004) enjoy sustained popularity because of their good performance and computational efficiency. Like ICA algorithms, the AJD approach allows extracting source components by groups, which appears to us an effective way to overcome the aforementioned limitation of assuming pair-wise spatial independence of all EEG source processes; source components may now be assumed independent between groups but not necessarily independent within each group. Mathematically, this amounts to require the products $\hat{B}Q_1\hat{B}^T$, $\hat{B}Q_2\hat{B}^T$, … be block-diagonal instead of diagonal. Such an approach has been foreseen by Cardoso (1998 b) and is nowadays referred to as *independent subspace analysis* (ISA). Block-AJD (B-AJD) allows seeking brain networks (groups of dependent source processes) instead of just several disjoint "hot spot", which is in line with current trends in brain neurophysiology (e.g., Mantini et al., 2007). As per today B-AJD is limited in practice by the necessity of specifying a-priori the numerosity and composition of the groups (Theis, 2005; Févotte and Theis, 2007). Research on AJD algorithms is currently flourishing. Recent trends include pursuing decomposition by blocks and seeking optimal weighting (e.g., Tichavskí et al., 2008). We believe that the resulting improvements hold promise for the BSS field and its applications to human electroencephalogram.

**SOS BSS methods solved by approximate joint diagonalization**

The first proposed SOS method (Féty and Uffelen, 1988; Tong et al., 1991 b) exploited signal coloration. It consisted on joint diagonalization of two matrices, the covariance matrix and a lagged covariance matrix, allowing an exact solution via the well-known generalized eigenvalue-eigenvector decomposition (Choi et al., 2002; Parra and Saida, 2003). The corresponding procedure for exploiting energy time variation traces back to the work of Souloumiac (1995); if the energy of a source component changes in two successive time intervals, then the component can be estimated by joint diagonalization of the two covariance matrices estimated on those intervals. Importantly, if the source is active in one interval and inactive in the other the obtained filter is optimal (Souloumiac, 1995). Along these lines see the discussion on *super-efficiency* in Pham and Cardoso (2001). Although very simple and fast, these two-matrix joint diagonalization methods are very sensitive to estimation errors of those matrices. If the noise covariance structure is different in the two matrices then the joint diagonalization of the signal structure is severely distorted. A considerable improvement is obtained by AJD of a larger set of matrices. This idea, first applied in SOBI for colored source components (Belouchrani et al., 1997), has been then applied to non stationary source components (Choi and Cichocki 2000; Pham and Cardoso, 2001) and finally extended to both colored and non stationary source components (Belouchrani and Amin, 1998; Pham, 2002).

Practically, in many SOS methods such as SOBI the data are first whitened and normalized (sometimes it is said they are sphered) as

$$z(t) = Hv(t),$$

where $H \in \mathbb{R}^{M \cdot N}$ is such that the covariance of $z(t)$ is the identity. Then, the AJD of a *set* of delayed covariance matrices (several lags: SOBI), and/or a set of covariance matrices on several windows of $z(t)$ is performed. It is known that the pre-whitening may jeopardize the separation performance due to the estimation error of the data covariance matrix, which is *exactly* diagonalized at the expense of the other matrices (Cardoso, 1994; Yeredor, 2000; Pham, 2001a). Hence, a better procedure is obtained by using robust whitening (Choi et al., 2002) or obtaining the AJD of a set of covariance matrices directly on $v(t)$, which amount to avoiding the pre-whitening step altogether (Ziehe and Müller, 1998), or by diagonalizing partial autocorrelation matrices (Dégerine and Malki, 2000). As compared to the two-matrix diagonalization the AJD approach is known to be more robust and efficient (Belouchrani et al., 1997; Belouchrani et Amin, 1998; Choi et al., 2002). One problem encountered by researchers with SOBI is how to choose an appropriate set of lags (Tang et al., 2004, 2005). For colored Gaussian auto-regressive (AR) processes the asymptotically optimal set of lags includes as many lags as necessary to describe the maximal order of the process (see for example Doron and Yeredor, 2004). The AR order/number of lags depending on several factors (e.g., sampling rate, number of peaks in the source power spectrum etc.), one should estimate it on data at hand. A simpler solution to this problem is treated in appendix (B). It arises after shifting the AJD problem into the time-frequency domain, framework that we now delineate.

**Time-frequency expansions**

Source separation methods can be applied in different representation spaces. In fact, applying to (1.0) any invertible and linearity-preserving transform $\mathcal{T}$ leads to

$$\mathcal{T}[v(t)] = A\mathcal{T}[s(t)],$$

which preserves the mixing model. Then, solving source separation in the transformed space still provides estimation of the matrix $A$ or of its inverse $B$, which can be used directly in Eq. (1.1) for recovering the source $s(t)$ in the initial space. For example, the transform $\mathcal{T}$ may be a discrete Fourier transform, a time-frequency transform such as the Wigner-Ville transform or a wavelet transform. AJD-based SOS methods such as SOBI can be easily and conveniently transposed in the frequency domain, thence in the time-frequency domain, whether we perform the frequency expansion for several time segments. Such approach is currently attracting much interest in the BSS community, especially for audio and speech applications (Belouchrani and Amin, 1998; Choi et al., 2002; Bousbia-Salah et al., 2003; Deville, 2003; Zhang and Amin, 2006; Aïssa-El-Bay et al., 2007). There exist several time-frequency expansions. For its simplicity in this study we consider the *short Fourier transform*, from which Fourier cospectral matrices are readily estimated[2]. We will compute Fourier cospectral matrices $\boldsymbol{C}_{(fi)} \in \mathbb{R}^{N \cdot N}$ for a range $f: 1…F$ of discrete frequencies and for a range $i: 1…I$ of temporal windows. Temporal windows should be short enough to capture the energy variations over time and wide enough to allow satisfactory estimations of cospectral matrices for each of them separately. For each temporal window $i$ the $f^{\text{th}}$ cospectral matrix $\boldsymbol{C}_{(fi)}$ holds the portion of the sensor covariance matrix corresponding to the $f^{\text{th}}$ frequency. Its diagonal elements hold the power (auto-spectra) of each measurement channel while its off-diagonal elements hold the terms describing the in-phase SOS *dependency* for that time window and frequency. As we have seen those off-diagonal terms are canceled by AJD in order to recover uncorrelated source components. Clearly, cospectral matrices are affine to the delayed covariance matrices used by SOBI, since they are a linear transformation of each other (e.g., Bloomfield, 2000, p. 12; Pham, 2001 a). Nonetheless, working in the frequency domain is advantageous for several reasons: first, covariance statistical estimations in the time domain are distorted for temporally correlated processes like EEG (Beran, 1994). Second, estimating cospectral matrices in the frequency domain is computationally more efficient than estimating delayed covariance matrices in the time domain[3]. Finally, the AJD of cospectra has been

---

[2] See appendix (A) for details on Fourier co-spectral matrices.
[3] We have analyzed the computational complexity of estimating the former and latter matrices. Fourier cospectra estimations may take advantage of efficient split-radix fast Fourier transform (FFT) algorithms such

connected to the *Gaussian mutual information* criterion (Pham 2001 a, 2002). This places the ensuing method at the hearth of the BSS theory and steers toward the Cramér-Rao bound (Pham, 2001 a; Pham and Cardoso, 2001). We are aware of only one study comparing the AJD of delayed covariance matrices (SOBI) to the AJD of cospectral matrices (Doron and Yeredor, 2004) and it clearly showed the better performance of the latter.

**Approximate joint diagonalization of cospectra (AJDC): an extended time-frequency approach**

Without loss of generality, the AJDC solution to the BSS problem (1.1) can be written compactly such as

$$\hat{\boldsymbol{B}} = AJD(\mathbb{C}), \qquad (1.6)$$

where $\mathbb{C} : \{\boldsymbol{C}_1, \boldsymbol{C}_2, \ldots\}$ is the *diagonalization set*, i.e., a set of estimated Fourier cospectral matrices to be simultaneously diagonalized. The rational behind AJDC is expressed schematically in Fig 1. Each cube of the parallelepiped in the figure represents abstractly a cospectral matrix. The grid of cubes represents the sampling of some source property unfolding along two continuous dimensions (time and frequency) and one discrete dimension (experimental conditions). The different pattern of shading in each cube represents the different *cospectral structure* of each sampled region of the defined space. If only one source component was involved the shading could be directly interpreted as color coded energy, but since in general many source components are considered we shall think at the shading as a coding for the dependency structure. The variations of the cospectral structure in the defined space along the dimensions are called *signatures*. We say that a source component has a *characteristic signature* if no other source component has the same signature. Successful separation of a source component is obtained if the diagonalization set describes a characteristic signature of it. In other words, the diagonalization set should include at least two matrices differing in the dependency structure of this source component (with respect to the others) and those changes must not be the same for any other source components. For example, if the source component is narrowband and its frequency range differs from the others (characteristic coloration), the cospectral structure of this

---

as FFTW3 (Frigo and Johnson, 2005); in typical situations we may expect the computation complexity of Fourier cospectral matrices be 20 to 100 times smaller as compared to lagged covariance matrices.

source component along the frequency dimension will change uniquely and this change will enable the identification of that component. Algorithms like SOBI seek those changes to recover source components having non proportional power spectrum. The advantage of the time-frequency approach is precisely that *either* coloration *or* non stationarity characteristic signature can be captured in the time-frequency plane and that either one suffices to achieve separation. Thus, the multidimensional approach is robust with respect to possible violations of each assumption taken separately. An important aspect of data expansion is that it enhances the characterization of source signatures; while the noise power tends to spread uniformly in the time-frequency plane the source power will concentrate in characteristic regions, thus the method is more robust with respect to noise as well (Belouchrani et Amin, 1998). The same arguments can be strengthened profiting of further source diversities simultaneously, such as those of physiological and experimental origin. For instance, to separate the posterior Alpha rhythms from the Rolandic Mu rhythms one may use the fact that posterior Alpha rhythms, but not Mu rhythms, are blocked by eyes opening (Niedermeyer, 2005 a). Two time intervals separated by eyes opening should then be considered. To exploit possible source energy diversity in several experimental conditions it suffices to average cospectral matrices separately within each condition, as indicated schematically in Fig. 1. This allows much flexibility, for an arbitrary number of cospectra computed on short time intervals can be averaged for each condition. To visualize the comprehensive nature of the method one may imagine the parallelepiped in Fig 1 in any number of dimensions susceptible to describe variations in some source statistical property.

Putting all this in mathematical formalism turns out simple and elegant. Without loss of generality we shall always proceed by (1.6) after defining

$$\mathcal{C} : \{ \boldsymbol{C}_{(\upsilon)} \}, \qquad (1.7)$$

where $\upsilon$ is just a container for an arbitrary number of indexes and where each index indicates the sampling along a dimension. For example, the diagonalization set of Fig. 1 is obtained by defining $\upsilon \triangleq fik$, where the cospectra at $F$ frequencies ($f : 1…F$) are estimated for $I$ time intervals ($i : 1…I$) and *K* experimental conditions (*k* : 1…*K*). With such a diagonalization set one would exploit *the diversity of source energy between conditions* in addition to generic coloration and time-varying energy; notice that in this case source components can be identified if their energy differs in at least two experimental conditions, *regardless* the uniqueness of their spectral and stationarity signatures (that is, even if the basic assumptions of the SOS BSS method do not hold), but also if their characteristic signature is across the frequency or time dimension but not across experimental conditions. Along the same line, we can exploit the reactivity of EEG oscillations as aforementioned discussed, the presence/absence of a steady-state sensory stimulation, the presence/absence of electrical or magnetic stimulation, etc.; one may add as many indexes as desired and always proceed by (1.6).

We have seen that adding dimensions for expanding the data increases the chance to uncover the characteristic signatures of source components and increase the robustness with respect to noise. However the number of matrices in the diagonalization set cannot be increased indefinitely. The essence of AJD algorithms consists in approximating the "average eigen-structure" of the input matrices. In general, any set (1.7) can be exactly jointly diagonalized if the instantaneous linear model holds exactly (Hyvärinen et al., 2001, p. 344), in which case all matrices in the set share common eigenvectors and the two subspaces spanned by those eigenvectors and the columns of mixing matrix *A* are identical (Belouchrani and Amin, 1998). This in practice will not quiet happen because of sampling estimation errors and noise, and while the latter is reduced the former is increased by data expansion, making more difficult finding the average eigen-structure. Another drawback of multiple dimension data expansion is that the instantaneous linear model (1.0) may not hold for all dimensions. Finally, an open question is how the time-frequency plane should be sampled (on the other hand sampling of experimental conditions is given by definition). We see that the proper definition of the diagonalization set is the very challenge of AJD-based algorithms. A useful tool to identify regions where the characteristic signatures reside is described in appendix (B).

**Discussion**

Blind source separation (BSS) is a widespread method used in a number of scientific and technical fields (Hyvärinen et al., 2001; Cichocki and Amari, 2002). Its use in EEG literature is currently growing at a fast pace. When applied to EEG data BSS decomposes scalp signals in a number of components. These components may correspond to the activity of cortical dipole layers generating the observed EEG. Precisely, BSS implicitly estimates their orientation and explicitly estimates their waveform (out of a sign and energy arbitrariness) and mixing coefficients. From the latter the spatial location can be estimated using an inverse solution method (Lopes da Silva, 2004; Greenblatt et al., 2005; see for example Van der Loo et al., 2007). Environmental and physiological artifacts may be extracted as well, while effective reduction of background noise may require additional noise suppression procedures (e.g., Vorobyov and Cichocki, 2002). Typically, BSS estimations feature higher SNR and strong suppression of the interference generated by other dipoles as compared to raw EEG. When the assumptions hold, BSS provides optimal estimations, in that they do not depend on physical modeling of the head. When the assumptions do not hold, BSS provides another representation of the sensor measurement space, which may still be useful (e.g., artifact reduction by BSS filtering as per Eq. 1.4), but that may encourage misleading interpretations of source waveform and associated topographies. One should be careful in claiming that the extracted components and associated topographies correspond to actual EEG physiological sources. Since checking the assumptions of the chosen BSS model and method is virtually impossible, the most credible arguments are those founded upon neurophysiological knowledge.

The linear BSS instantaneous model (1.0) makes a number of restrictive assumptions that are rarely checked or investigated. One assumption is that the number of sources is not greater than the number of sensors. When this is not the case (undetermined case) it is not possible to solve the BSS problem unless other constrains on the sources are introduced (e.g., sparsity: Gribonval and Lesage, 2006). One may also wonder if during the analyzed time interval the number of active dipoles is stable (Li et al., 2006). In practice, brain electrical "source components" are macroscopic electric dipole with relatively high SNR formed by the synchronous activity of pyramidal cells over large cortical areas (Nunez and Srinivasan, 2006). For sufficiently small time intervals one may assume that such high-SNR layers are limited in number. Other concurrently active cortical columns may be ignored if their current is comparatively negligible and it does not matter if the dipoles are active throughout the time interval or intermittently (actually such non stationarity signature can be explicitly exploited). Henceforth, assuming at least as many sensors as relevant sources does not appear problematic if we consider a sufficiently small time interval. Still, no definitive solution exists to the problem of estimating the number of source components in the overdetermined case (more sensor than source components). Whereas correct dimensionality reduction (appendix E) allows exact determination, over-reduction must be avoided since in this case identifiability is lost and several generators are extracted mixed in one component. A safe strategy is to identify a few meaningful components and keep reducing the dimension until those components are not distorted.

Another assumption is that the mixing matrix $A$ in (1.0) is full-column rank. The columns of $A$ are scalp spatial pattern vectors of the source components and the more the electrodes are close to each other, the more those vectors will be collinear. Consequently, it is always better to space the electrodes as much as possible on the scalp[4]. Several restrictive assumptions are made by model (1.0) also on the nature of brain electric fields. One may ask whether it is reasonable to assume that dipoles keep fixed orientation and location in the analyzed time interval. Each row vector of the matrix $B$ can be conceived as a *spatial filter* extracting the electric field of a dipole with a given fixed spatial extension, location and orientation. For a fixed spatial sensor configuration with respect to the brain, which is the case of a single EEG recording session, the orientation and location of electric dipoles are fixed by the anatomy and physiology of the grey matter forming the dipole. However, the dipole approximation becomes untenable for sources distributed over large areas (Malmivuo and Plonsey, 1995; Nunez and Srinivasan 2006). Also, there is convincing evidence of traveling waves phenomena in the brain; long wavelength waves originating in a region and propagating via cortico-cortical connections to other regions (Lopes da Silva and Van Rotterdam, 2005; Srinivasan et al., 2006;

---

[4] This suggests that placing many electrodes closely spaced above the brain region of interest, as it is sometimes done, is not a convenient strategy if multivariate statistical methods are to be employed.

Thorpe et al., 2007). These phenomena cannot be modeled by an instantaneous model and become more equivocal with larger time intervals. Also, the longer the time interval under analysis the less tenable is the stationarity assumption, which is basic to SOS estimations (Hyvärinen et al., 2001, p. 49). At the same time one must care to retain enough data points for analysis in order to avoid *overfitting* (Müller et al., 2004). Särelä and Vigário (2003) reported that using small time intervals the output may contain artefacts that are not present in the data. For HOS method such as FastICA artifacts takes the form of artificial spikes and bumps, whereas for SOS methods such as SOBI they take the form of artificial sinusoid waves. Meinecke et al. (2002) and Müller et al. (2004) addressed the problem of obtaining robust and reliable source estimates. They proposed a resampling-based methods consisting in running the algorithms on different time intervals and retain only the source processes that can be found consistently. In conclusion, although statistical estimations improve with the number of samples we advocate the use of multiple time intervals as short as possible (enough to avoid overfitting while justifying the BSS method assumptions), modeling appropriately the stationarity within intervals while exploiting explicitly the non stationarity between intervals. In this sense an efficient time-frequency approach appears a precious option. Although we have contended that SOS-AJD methods such as AJDC fit well EEG data in general, a safe strategy is to compare the output to at least one HOS methods with any real EEG data problem at hand. We also notice that in the EEG field the instantaneous model has been rarely challenged (Anemüller et al., 2003; Dyrholm et al., 2005). It is unfortunate that throughout comparisons of linear instantaneous, time-varying and convolutive model are lacking since the latter two families of BSS models may admit moving dipoles and traveling waves.

In this study we have described a simple time-frequency approach based on the approximate joint diagonalization of Fourier cospectral matrices (AJDC). AJDC is an extension of popular AJD-based algorithms such as SOBI, which are derived as restricted instances (exploiting source coloration or source non stationarity only), yet it is efficient statistically and computationally. Computationally, the AJDC equivalent of SOBI is several tens of times faster than SOBI. In turn, SOBI is known to be faster than, in the order, JADE, FastICA and InfoMax, the latter being the slowest (Kachenoura et al., 2008). Although those authors do not quantify precisely the complexity of each algorithm, we can safely say that AJDC is tens of times faster than SOBI and hundreds to thousands times faster than FastICA and InfoMax. Moreover AJDC (as all AJD-based algorithms) does not require parameter tuning for convergence. However, it requires an appropriate definition of the diagonalization set to correctly identifying the potential diversities in the data set. Instead of understanding this as a nuisance, we have contended that it amounts to correctly identifying the relevant aspects of the data variance at hand. Such an "informed" approach is somehow in between the completely blind setting, in which no a-priori knowledge on the source is assumed and the semi-blind approach, where temporal, spatial, spectral or other constraints are introduced in the cost function (Roberts, 1998; Ille, Berg and Scherg 2002; James and Gibson, 2003, James and Hesse, 2005; Lu and Rajapakse, 2005; Hesse and James, 2006; Barbati et al., 2006; Wang and James, 2007; Barbati et al, 2008; Zhang, 2008). The basic time-frequency approach exploits the temporal dependency and energy variation over time of EEG. The diagonalization scheme can be defined so as to maximize the chance of separating dipole layers responsible for brain functions studied by experimental manipulation. Assessing the difference in two or more experimental conditions is customary in cognitive and clinical studies using either continuous recording or evoked potentials paradigms. In this sense, AJDC may be an ideal companion for a very wide range of EEG experimental research. In the appendix we have collected several useful details about the effective use of AJDC for EEG data, which cannot be found elsewhere. Those details may be valuable to the reader interested in implementing AJDC or other time-frequency BSS algorithms. Code for AJD algorithms is commonly publicly available. An executable application performing BSS by AJDC can be obtained upon request to the corresponding author.

**References**

Aissa-El-Bey A, Linh-Trung N, Abed-Meraim K, Belouchrani A, Grenier Y. Underdetermined Blind Separation of Nondisjoint Sources in the Time-Frequency Domain. IEEE Trans Signal Process 2007; 55(3): 897-907.

Anemüller J, Sejnowski TJ, Makeig S. Complex independent component analysis of frequency-domain electroencephalographic data. Neural Netw 2003; 16(9): 1311-23.

Ans B, Hérault J. Jutten C. Adaptive Neural Architectures: Detection of Primitives. In : Proc. COGNITIVA 1985 : 593-597.

Barbati G, Sigismondi R, Zappasodi F, Porcaro C, Graziadio S, Valente G, et al. Functional source separation from magnetoencephalographic signals. Hum Brain Mapp 2006; 27(12): 925-934.

Barbati G, Porcaro C, Hadjipapas A, Adjamian P, Pizzella V, Romani GL, et al. Functional source separation applied to induced visual gamma activity. Hum Brain Mapp 2008; 29(2): 131-141.

Bell AJ, Sejnowski TJ. An Information-Maximization Approach to Blind Separation and Blind Deconvolution. Neural Comput 1995; 7: 1129-1159.

Belouchrani A, Abed-Meraim K, Cardoso J-F, Moulines E. A blind source separation technique using second-order statistics. IEEE Trans Signal Process 1997; 45(2): 434-444.

Belouchrani A, Amin MG. Blind Source Separation Based on Time-Frequency Signal Representations. IEEE Trans Signal Process 1998; 46(11): 2888-2897.

Beran J. Statistics for Long-Memory processes. London: Chapman & Hall, 1994.

Bloomfield P. Fourier Analysis of Time Series. New York: John Wiley & Sons, 2000.

Bousbia-Salah A, Belouchrani A, Bousbia-Salah H. A one step time-frequency blind identification. 7[th] Int Symp Sig Process Applications 2003; 1(1): 581- 584.

Buzsáki G. Rhythms of the Brain. New York: Oxford Univ Press, 2006.

Cao J, Murata N, Amari S-I, Cichocki A, Takeda T. Independent component analysis for anaveraged single-trial MEG data decomposition and single-dipole localization. Neurocomputing 2002; 49: 255-277.

Cardoso JF, Souloumiac A. Blind beamforming for non-Gaussian signals. IEE Proc-F (Radar and Signal Process) 1993; 140(6): 362-370.

Cardoso J.-F. On the performance of orthogonal source separation algorithms. Proc EUSIPCO, Edinburg (UK) 1994; 776-779.

Cardoso J-F. Blind Signal Separation: Statistical Principles. IEEE Proc 1998a; 9(10): 2009-2025.

Cardoso J-F. Multidimensional independent component analysis. Proc ICASSP 1998b, Seattle, USA: 1941-1944.

Cardoso J-F. High-Order Contrasts for Independent Component Analysis. Neural Comput 1999; 11(1): 157-192.

Choi S, Cichocki A. Blind Separation of nonstationary sources in noisy mixtures. Electron Lett 2000; 36: 848-849.

Choi S, Cichocki A, Belouchrani. Second Order Nonstationary Source Separation. J VLSI Sig. Process 2002; 32(1-2): 93-104.

Choi S, Cichocki A, Park HM, Lee S-Y. Blind source separation and independent component analysis: A review, Neural Inf Process - Letters and Reviews 2005, 6(1): 1-57.

Cichocki A, Amari SI. Adaptive Blind Signal and Image Processing. Learning Algorithms and Applicaions. New-York: John Wiley & Sons, 2002.

Congedo M. Subspace Projection Filters for Real-Time Brain Electromagnetic Imaging. IEEE Trans Biomed Eng 2006; 53(8): 1624-34.

Cotter S, Rao BD, Engan K, Kreutz-Delgado K. Sparse solutions to linear inverse problems with multiple measurement vectors. IEEE Trans Signal Process 2005; 53(7): 2477-2488.

Crespo-Garcia M, Atienza M, Cantero JL. Muscle artifact removal from human sleep EEG by using independent component analysis. Ann Biomed Eng 2008; 36(3): 467-75.

Darmois G. Analyse générale des liaisons stochastiques. Rev Inst Inter Stat 1953; 21, 2-8.

Dat TH, Guan C. Feature Selection Based on Fisher Ratio and Mutual Information Analyses for Robust Brain Computer Interface. Proc ICASSP 2007; 1: I-337-I-340.

De Clercq W, Vergult A, Vanrumste B, Van Paesschen W, Van Huffel S. Canonical correlation analysis applied to remove muscle artifacts from the electroencephalogram. IEEE Trans Biomed Eng 2006; 53(12 Pt 1):2583-2587.

Dégerine S, Malki R. Second-Order Blind Separation of Sources Based on Canonical partial Innovations. IEEE Trans Signal Process 2000; 48(3): 629-641.

Degerine S, Kane E. A Comparative Study of Approximate Joint Diagonalization Algorithms for Blind Source Separation in Presence of Additive Noise. IEEE Trans Signal Process 2007; 55 (6-2): 3022-3031.

Deville Y. Temporal and time-frequency correlation-based blind source separation methods, Proc ICA2003, Nara, Japan, 2003: 1059-1064.

Doron E, Yeredor A. Asymptotically Optimal Blind Separation of Parametric Gaussian Sources. Proc ICA 2004, C.G Puntonet and A. Prieto (Eds.): 390-397.

Dyrholm M, Makeig S, Hansen LK. Model selection for convolutive ICA with an application to spatiotemporal analysis of EEG. Neural Comput 2007; 19(4): 934-55.

Fadaili EM, Moreau NT, Moreau E. Nonorthogonal Joint Diagonalization/Zero Diagonalization for Source Separation Based on Time-Frequency Distributions. IEEE Trans Signal Process 2007; 55(5-1): 1673-1687.

Féty L, Uffelen J.-P. New Methods for Signal Separation. Proc. of the 14[th] conf. on HF Radio System and Techniques, London, 1988; 226-230.

Févotte C, Theis FJ. Pivot selection strategies in jacobi joint block-diagonalization. In Proc. ICA, London, U.K 2007; 4666 of LNCS: 177-184.

Fitzgibbon SP, Powers DM, Pope KJ, Clark CR. Removal of EEG noise and artifact using blind source separation. J Clin Neurophysiol 2007; 24(3): 232-243.

Frank RM, Frishkoff GA. Automated protocol for evaluation of electromagnetic component separation (APECS): Application of a framework for evaluating statistical methods of blink extraction from multichannel EEG. Clin Neurophysiol 2007; 118(1): 80-97.

Frigo M, Johnson SG. The Design and Implementation of FFTW3. Proc IEEE 2005; 93(2), 216-231.

Gasser T, Schuller JC, Gasser US. Correction of muscle artefacts in the EEG power spectrum. Clin Neurophysiol 2005; 116(9): 2044-50.

Gorodnitsky IF, George JS, Rao BD. Neuromagnetic source imaging with FOCUSS: a recursive weighted minimum norm algorithm. Electroencephalogr Clin Neurophysiol 1995; 95(4):231-51

Greenblatt RE, Ossadtchi A, Pflieger ME. Local Linear Estimators for the Bioelectromagnetic Inverse Problem. IEEE Trans Signal Process 2005; 53(9): 3403-3412.

Gribonval R, Lesage S. A survey of Sparse Component Analysis for blind source separation: principles, perspectives, and new challenges. Proc. of Eur Symp Artif Neural Netw (ESANN 2006) 2006: 323-330.

Guimaraes MP, Wong DK, Uy ET, Grosenick L, Suppes P. Single-trial classification of MEG recordings. IEEE Trans Biomed Eng 2007; 54(3): 436-443.

Halder S, Bensch M, Mellinger J, Bogdan M, Kübler A, Birbaumer N, et al. Online artifact removal for brain-computer interfaces using support vector machines and blind source separation. Comput Intell Neurosci 2007: 82069.

Hérault J., Jutten C. Space or time adaptive signal processing by neural network models. Proc Int Conf Neural Netw Computing, Snowbird (Utah), April 1986; 151, 206-211.

Hesse CW, James CJ. On semi-blind source separation using spatial constraints with applications in EEG analysis. IEEE Trans Biomed Eng 2006; 53(12): 2525-2534.

Hyvärinen A. Fast and robust fixed-point algorithms for independent component analysis. IEEE Trans Neural Netw 1999; 10(3): 626-634.

Hyvärinen A, Karhunen J, Oja E. Independent Component Analysis. New York: John Wiley & Sons, 2001.

Ille N, Berg P, Scherg M. Artifact correction of the ongoing EEG using spatial filters based on artifact and brain signal topographies. J Clin Neurophysiol 2002; 19(2): 113-124.

Iriarte J, Urrestarazu E, Valencia M, Alegre M, Malanda A, Viteri C, Artieda J. Independent component analysis as a tool to eliminate artifacts in EEG: a quantitative study. J Clin Neurophysiol 2003; 20(4): 249-257.

James C, Gibson O. Temporally constrained ICA: an application to artifact rejection in electromagnetic brain signal analysis. IEEE Trans Biomed Eng 2003; 50(9): 1108-1116.

James C, Hesse C. On the use of Spectrally Constrained ICA applied to single-channel Ictal EEG recordings within a Dynamical Embedding Framework. Conf Proc IEEE Eng Med Biol Soc 2005; 1: 956-659.

Joyce CA, Gorodnitsky IF, Kutas M. Automatic removal of eye movement and blink artifacts from EEG data using blind component separation. Psychophysiology 2004; 41(2): 313-25.

Jung TP, Makeig S, Westerfield M, Townsend J, Courchesne E, Sejnowski TJ. Removal of eye activity artifacts from visual event-related potentials in normal and clinical subjects. Clin Neurophysiol 2000; 111(10): 1745-58.

Jutten C, Herault J. Blind separation of sources, Part 1: an adaptive algorithm based on neuromimetic architecture. Signal Process 1991; 24(1): 1-10.

Kachenoura A, Albera L, Senhadji L, Comon P. ICA: a potential tool for BCI systems. IEEE signal Process Mag 2008; 25(1):57-68.

Kierkels JJM, van Boxtel GJM, Vogten LLM. A model-based objective evaluation of eye movement correction in EEG recordings. IEEE Trans Biomed Eng 2006; 53(2): 246-253.

Lee DD, Seung HS. Learning the Parts of Objects by Nonnegative Matrix Factorization, Nature 1999; 401: 788-791.

Lee DD, Seung HS. Algorithms for non-negative matrix factorization. Adv Neural Info Proc Syst 2001; 13: 556-562.

Lemm S, Curio G, Hlushchuk Y, Müller K-R. Enhancing the Signal-to-Noise Ratio of ICA-Based Extracted ERPs. IEEE Trans Biomed Eng 2006, 53(4): 601-607.

Li X-L, Zhang X-D. Nonorthogonal Joint Diagonalization Free of Degenerate Solution. IEEE Trans Sig Proc 2007; 55(5): 1803-1814.

Li Y, Cichocki A, Amari S-I. Blind estimation of channel parameters and source components for EEG signals: a sparse factorization approach. IEEE Trans Neural Netw 2006; 17(2): 419-431.

Lopes da Silva F. Functional Localization of Brain Sources using EEG and/or MEG data: Volume Conductor and Source Models. Magn Res Img 2004, 22: 1533-1538.

Lopes da Silva F. Computer-Assisted EEG diagnosis: Pattern Recognition and Brain Mapping. In: Electroencephalography. Basic Principles, Clinical Applications, and Related Fields. Niedermeyer E and Lopes da Silva F. (Eds), 5th ed., New York: Lippincott Williams & Wilkins, 2005 a: 1233-1263.

Lopes da Silva F. Event Related Potentials: Methodology and Quantification. In: Electroencephalography. Basic Principles, Clinical Applications, and Related Fields. Niedermeyer E and Lopes da Silva F. (Eds), 5th ed., New York: Lippincott Williams & Wilkins, 2005 b: 991-1001.

Lopes da Silva F. Van Rotterdam A. Biophysical Aspects of EEG and Magnetoencephalogram Generation. In: Electroencephalography. Basic Principles, Clinical Applications, and Related Fields. Niedermeyer E and Lopes da Silva F. (Eds), 5th ed., New York: Lippincott Williams & Wilkins, 2005: 107-125.

Lu W, Rajapaske JC. Approach and Applications of Constrained ICA. IEEE Trans Neural Netw 2005; 16(1): 203-212.

Matsuoka K., Ohya M., Kawamoto M. A neural net for blind separation of nonstationary signals. Neural Netw 1995; 8:3, 411-419.

Malmivuo J, Plonsey R. Bioelectromagnetism. Principles and Applications of Bioelectric and Biomagnetic Fields. New York: Oxford Univ Press, New York, 1995.

Mantini D, Perrucci, MG, Del Gratta C, Romani, GL, Corbetta M.Electrophisiological signatures of resting state networks in the human brain. PNAS 2007; 104(32):10.1073.

Meinecke F, Ziehe A, Kawanabe M, Müller K R. A resampling approach to estimate the stability of one- or multidimensional independent components. IEEE Trans Biomed Eng 2002; 49, 1514-1525.

Molgedey L, Schuster, HG Separation of a Mixture of Independent Signals using Time Delayed Correlations. Phys Rev Lett 1994; 72: 3634-3636.

Mosher JC, Lewis PS, Leahy RM. Multiple dipole modeling and localization from spatio-temporal MEGdata. IEEE Trans Biomed Eng 1992; 39(6): 541-557.

Müller KR, Vigario R, Meinecke F, and Ziehe A. Blind source separation techniques for decomposing event-related brain signals. Int J Bifurcat Chaos 2004; 14(2): 773-791.

Niedermeyer E. The Normal EEG of the waking Adult. In: Electroencephalography. Basic Principles, Clinical Applications, and Related Fields. Niedermeyer E and Lopes da Silva F. (Eds), 5th ed., New York: Lippincott Williams & Wilkins, 2005 a: 167-191.

Niedermeyer E. Sleep and EEG. In: Electroencephalography. Basic Principles, Clinical Applications, and Related Fields. Niedermeyer E and Lopes da Silva F. (Eds), 5th ed., New York: Lippincott Williams & Wilkins, 2005 b: 193-207.

Niedermeyer E. Epileptic Seizure Disorders. In: Electroencephalography. Basic Principles, Clinical Applications, and Related Fields. Niedermeyer E and Lopes da Silva F. (Eds), 5th ed., New York: Lippincott Williams & Wilkins, 2005 c: 505-619.

Nunez PL, Srinivasan R. Electric Field of the Brain, 2nd ed., New York: Oxford Univ Press, 2006.

Parra L, Sajda P. Blind Source Separation via Generalized Eigenvalue Decomposition, J Mach Learn Res 2003; 4: 1261-1269.

Pham D-T. Blind Separation of Instantaneous Mixture of Sources via the Gaussian Mutual Information Criterion. Signal Process 2001a; 81: 855-870.

Pham D-T; Joint Approximate Diagonalization of Positive Definite Matrices. SIAM J. on Matrix Anal and Appl 2001b; 22(4): 1136-1152.

Pham D-T, Cardoso J-F. Blind Separation of Instantaneous Mixtures of Non Stationary Sources. IEEE Trans Signal Process 2001; 49(9): 1837-1848.

Pham D-T. Exploiting source non stationary and coloration in blind source separation. Digital Signal Process 2002; 1: 151-154.

Phlypo R, Boon P, D'Asseler Y, Lemahieu I. Removing ocular movement artefacts by a joint smoothened subspace estimator. Comput Intell Neurosci 2007: ID 75079.

Pfurtscheller G, Lopes da Silva F.H. Event-related EEG/MEG synchronization and desynchronization: basic principles. Clin Neurophysiol 1999; 110(11): 1842-57.

Qin L, Ding L, He B. Motor Imagery Classification by Means of Source Analysis for Brain-Computer Interface Applications. J Neural Eng 2004; 1: 135-141.

Roberts S. Independent component analysis: source assessment and separation, a Bayesian approach. IEE Proc Vis Image Signal Process 1998; 145(3): 149-154.

Rodríguez-Rivera A, Baryshnikov BV, Van Veen BD, Wakai RT. MEG and EEG Source Localization in Beamspace. IEEE Trans Biomed Eng 2006; 53(3): 430-441.

Romero S, Mañanas MA, Barbanoj MJ. A comparative study of automatic techniques for ocular artifact reduction in spontaneous EEG signals based on clinical target variables: A simulation case. Comput Biol Med. 2008; 38(3): 348-360.

Sander TH, Burghoff M, Curio G, Trahms L. Single Evoked Somatosensory MEG Responses Extracted by Time Delayed Decorrelation. IEEE Trans Signal process 2005; 53(9): 3384-3392.

Särelä J, Vigário R. Overlearning in Marginal Distribution-Based ICA: Analysis and Solutions. JMach Learn Res 2003; 4: 1447-1469.

Sarvas J. Basic Mathematical and Electromagnetic Concepts of the Biomagnetic Inverse Problem. Phys Med Biol 1987; 32(1): 11-22.

Serby H. Yom-Tov E. Inbar GF. An improved P300-Brain-Computer Interface. IEEE Trans Neural Syst Rehabil Eng 2005; 13(1): 89-98.

Souloumiac A. Blind Source Detection and separation using second order nonstationarity. In Proc ICASSP 1995; 1912-1915.

Speckmann E-J, Elger CE. Introduction to the Neurophysiologicalal Basis of the EEG and DC Potentials. In: Electroencephalography. Basic Principles, Clinical Applications, and Related Fields. Niedermeyer E and Lopes da Silva F. (Eds), 5th ed., New York: Lippincott Williams & Wilkins, 2005: 17-29.

Srinivasan R, Winter WR, Nunez PL. Source analysis of EEG oscillations using high-resolution EEG and MEG; Prog Brain Re. 2006; 159: 29–42.

Steriade M. Cellular Substrates of Brain Rhythms. In: Electroencephalography. Basic Principles, Clinical Applications, and Related Fields. Niedermeyer E and Lopes da Silva F. (Eds), 5th ed., New York: Lippincott Williams & Wilkins, 2005: 31-83.

Sutherland MT, Tang. AC. Reliable detection of bilateral activation in human primary somatosensory cortex by unilateral median nerve stimulation. Neuroimage 2006, 33: 1024-1054.

Tang AC, Sutherland MT, McKinney CJ. Validation of SOBI components from High-density EEG. Neuroimage 2004, 25: 539-553.

Tang AC, Liu J-Y, Sutherland MT. Recovery of correlated neuronal sources from EEG: The good and bad ways of using SOBI. Neuroimage 2005, 28: 507-519.

Tang AC, Sutherland MT, Wang Y. Contrasting single-trial ERPs between experimental manipulations: Improving differentiability by blind source separation. Neuroimage 2006, 29: 335-346.

Theis FJ. Uniqueness of Complex and Multidimensional Independent Component Analysis. Signal Process 2004; 84: 951-956.

Theis FJ. Blind signal separation into groups of dependent signals using joint block diagonalization. In Proc. ISCAS, Kobe, Japan, 2005: 5878-5881.

Theis FJ, Inouye Y. On the use of joint diagonalization in blind signal processing. In Proc. ISCAS, Kos, Greece, 2006.

Thorpe SG, Nunez PL, Srinivasan R. Identification of wave-like spatial structure in the SSVEP: Comparison of simultaneous EEG and MEG. Stat Med 2007; 26: 3911-3926.

Tichavskí P, Yeredor A, Nielsen J. A Fast Approximate Joint Diagonalization Algorithm using a Criterion with a Block Diagonal Matrix. Proc ICASSP 2008, Las Vegas, USA, April 2008.

Tong L, Liu RW, Huang, Y-F Blind Estimation of correlated source signals. Sig Syst Computers 1990; 1: 258-262.

Tong L, Liu RW, Soon, VC, Huang Y-F. Indeterminacy and Identifiability of Blind Identification. IEEE Trans Circuits Syst 1991 a; 38(5): 499-509.

Tong L, Soon V, Huang Y. Liu RW. A necessary and sufficient condition Waveform-Preserving Blind Estimation of Multiple Independent Sources. IEEE Trans Signal Process 1991 b; 41(7): 2461-2470.

Tong L, Inouye Y, Liu RW. Waveform-Preserving Blind Estimation of Multiple Independent Sources IEEE Trans Signal Process 1993; 41(7): 2461-2470.

Van Der Loo E, Congedo M, Plazier M, Van De Heyning P, De Ridder D. Correlation between Independent Components of scalp EEG and intra-cranial EEG (iEEG) time series Int J Bioelectromagnetism 2007; 9, 4: 270-275.

Vigário RN. Extraction of ocular artifacts from EEG using independent component analysis. Electroenceph Clin Neurophysiol 1997; 103: 395-404.

Vorobyov S, Cichocki A. Blind Noise Reduction for Multisensory Signals using ICA and Subspace Filtering, with Applications to EEG Analysis. Biol Cybern 2002; 86: 293-303.

Vollgraf R, Obermayer K. Quadratic Optimization for Simultaneous Matrix Diagonalization. IEEE Trans Sig Process 2006; 54(9): 3270- 3278.

Wang S, James CJ. Extracting Rhythmic Brain Activity for Brain-Computer Interfacing through Constrained Independent Component Analysis. Comput Intell Neurosci 2007; ID 41468.

Wax M, Sheinvald J. A Least-Squares Approach to Joint Diagonalization. IEEE Signal Process Lett 1997; 4(2) 52-53.

Welch PD. The Use of Fast Fourier Transform for the Estimaton of Power Spectra: A Method Based on Time Averaging Over Short, Modified Periodograms. IEEE Trans Audio Electroacoustics 1967; 15(2): 70-74.

Whitham EM, Pope KJ, Fitzgibbon SP, Lewis T, Clark CR, Loveless S, et al. Scalp electrical recording during paralysis: quantitative evidence that EEG frequencies above 20 Hz are contaminated by EMG. Clin Neurophysiol 2007; 118(8): 1877-88.

Yeredor A. Blind Separation of Gaussian Sources via Second-Order Statistics with Asymptotically Optimal Weighting, IEEE Signal Process Lett 2000; 7(7): 197-200.

Yeredor A. Non-orthogonal joint diagonalization in the least-squares sensewith application in blind source separation. . IEEE Trans Signal Process 2002; 50 (7): 1545-1553.

Zhang Z-L. Morphologically constrained ICA for extracting weak temporally correlated signals. Neurocomputing 2008; 71: 1669-1679.

Zhang Y and Amin M. Blind separation of nonstationary sources based on spatial time-frequency distributions. EURASIP J Appl Signal Process 2006; ID 64785.

Zhou W, Gotman J. Removal of EMG and ECG artifacts from EEG based on wavelet transform and ICA. Conf Proc IEEE Eng Med Biol Soc 2004; 1: 392-5.

Zeman PM, Till BC, Livingston NJ, Tanaka JW, Driessen PF. Independent component analysis and clustering improve signal-to-noise ratio for statistical analysis of event-related potentials. Clin Neurophysiol 2007; 118(12): 2591-2604.

Ziehe A, Müller K-R. TDSEP–an efficient algorithm for blind separation using time structure. Proc Int Conf Artif Neural Netw (ICANN'98) 1998: 675–680.

Ziehe A, Laskov P, Nolte G, Müller R-K. A Fast Algorithm for Joint Diagonalization with Non Orthogonal Transformations and its Application to Blind Source Separation. J Mach Learn Res 2004; 5: 777-800.

BSS of Human EEG by SOS AJD– Congedo et al. 2008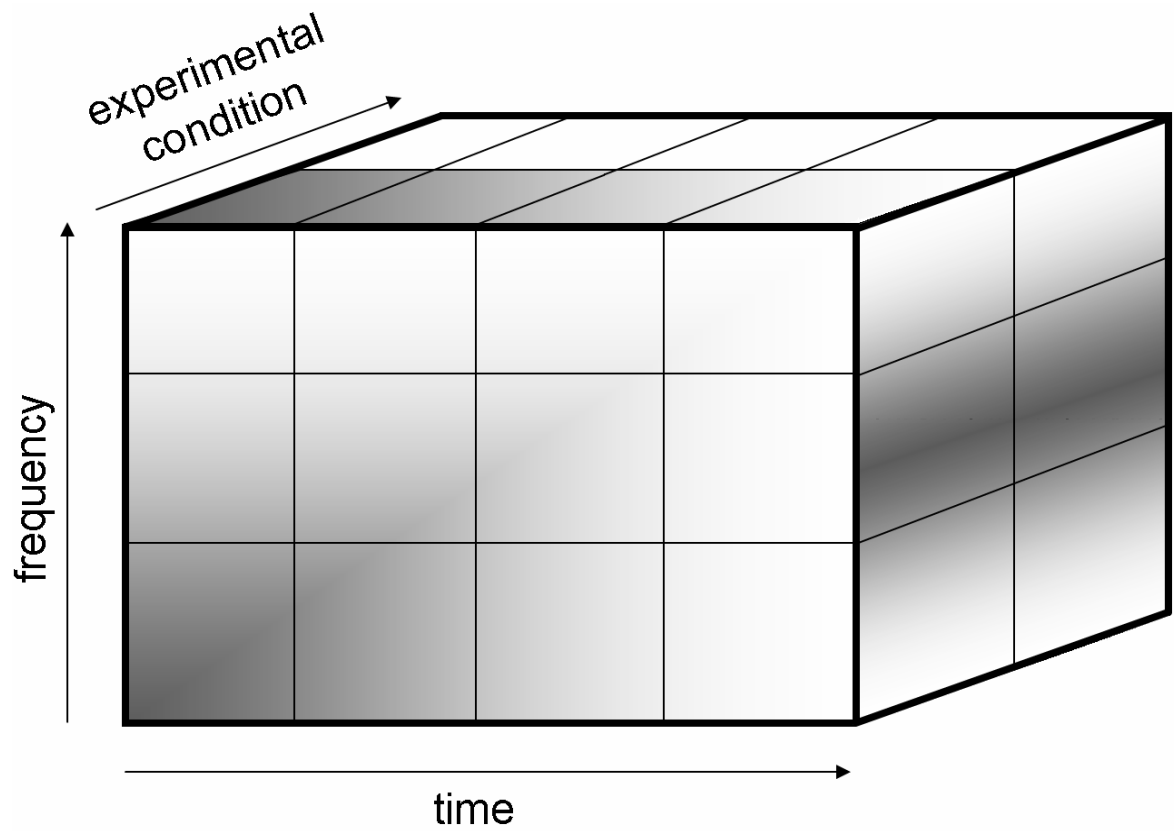

Figure 1. Schematic representation of the extended time-frequency blind source separation approach based on approximate joint diagonalization. See text for details.

**Appendix**

*A. Fourier cospectral matrices*

For real data, if *L* is the length of the time window and *S* the sampling rate expressed in Hz, there are $F=L/2+1$ Fourier frequencies with resolution $r=S/L$ equally spacing the range from 0Hz (DC-level) to the folding frequency (0Hz, 1*r* Hz, 2*r* Hz,…,*L/2r* Hz). Typically, we take both *S* and *L* as a power of two and cospectral estimates may be averaged within arbitrary time intervals by sliding overlapping windows. The latter strategy allows arbitrary time intervals length. The Fourier *cospectra* and *quadrature spectra* are defined as the real and imaginary part of the Fourier *cross-spectra* (Bloomfield, 2000). They are estimations of, respectively, the in-phase (or with a half cycle phase shift, i.e., opposite sign) and out-of-phase (a quarter cycle in either direction) covariance structure at frequency *f*. The discrete Fourier transform of sampled time-series $x_{(t)}$ over an epoch of length *L* is given by

$$d_x(f) = \frac{1}{L}\sum_{t=0}^{L-1} x(t) e^{-2\pi i f t}.$$

Let $d_x^{\Re}(f)$ and $d_x^{\Im}(f)$ be the real and imaginary part of $d_x(f)$, respectively[5]. Those coefficients are readily and efficiently estimated by fast Fourier Transform (FFT: Cooley and Tukey, 1965; Frigo and Johnson, 2005). Here below is the formula for computing the *cospectral matrix* at frequency *f* for time-series $x_{(t)}$ and $y_{(t)}$:

$$\boldsymbol{C}_{(f)} = \begin{pmatrix} d_x^{\Re}(f)d_x^{\Re}(f) + d_x^{\Im}(f)d_x^{\Im}(f) & d_x^{\Re}(f)d_y^{\Re}(f) + d_x^{\Im}(f)d_y^{\Im}(f) \\ d_y^{\Re}(f)d_x^{\Re}(f) + d_y^{\Im}(f)d_x^{\Im}(f) & d_y^{\Re}(f)d_y^{\Re}(f) + d_y^{\Im}(f)d_y^{\Im}(f) \end{pmatrix}.$$

The formula readily extends to any *N*-dimensional input time-series to obtain its *N*-dimensional cospectral matrix. Notice that the cospectral matrix is symmetric and that the diagonal elements are the auto-spectra, better knows as power spectra. For an arbitrary long EEG segment we typically obtain an estimate of the cospectral matrix averaging $\boldsymbol{C}_{(f)}$ over overlapping epochs of length *L* (Welch, 1967). Such estimates may then be summed across adjacent frequencies to obtain estimates within band-pass regions of interest. Summing all of them yields the covariance matrix as per Parseval's theorem. Pham (2001 a) shows that its Gaussian mutual information approximate joint diagonalization (AJD) criterion is the same using (complex) cross-spectral matrices or their real part (cospectral matrices). More in general, estimation of the separating matrix **B** can be obtained by AJD of cospectral matrix, quadrature spectra matrices or cross-spectral matrices if the mixing matrix is real (Theis, 2004), which we always assume as per (1.0) due to the absence of capacitive effects in the brain (see introduction).

---

[5] For the first (0Hz) and last ($L/2\kappa$Hz) Fourier frequency the coefficients are real.

*B. Defining diagonalization sets*

Some AJD algorithms allow weighting the diagonalization effort across the input matrices (e.g., Pham, 2001 b). Here we seek an adaptive solution to the weighting problem. Let us consider a set **w** of non-negative weights to be associated to each matrix of the diagonalization set. For consistency, we take **w** so that the average of its elements equals 1.0. As a general strategy, we may encourage the mask to be sparse (many zero entries), so to enable data expansion in multiple dimensions while keeping the diagonalization set of reasonable size. Because of volume conduction dipole fields result in covariance structures with many non-null off-diagonal terms (see for example simulated forward solutions in Congedo, 2006). On the other hand, higher frequency cospectra tend to be near diagonal because EEG energy decreases with frequency while spatially uncorrelated noise features a diagonal covariance structure. Whitham et al. (2007) even suggest that scalp EEG recording above 20Hz contains mainly electromyographic (EMG) activity. EMG exhibits near-diagonal covariance structures because it is spatially focal and does not propagate easily (through skull) to other leads. Consequently, matrices close to diagonal form should be down-weighted. In general, down-weighting low signal-to-noise ratio matrices amounts to effective noise suppression for the estimation of the separating matrix (Pham, 2001 a). Therefore, let us define for any cospectra $C_{(v)}$

$$\delta(C_{(v)}) = \frac{1}{N-1} \frac{\sum_{r \neq c} C_{(v)rc}^2}{\sum_{r=c} C_{(v)rc}^2}, \quad (B.1)$$

where $C_{(v)rc}$ is the entry of matrix $C_{(v)}$ at row $r$ and column $c$ and $N$ is the size of the matrix (number of channels). For a positive definite matrix, measure (B.1) is bounded inferiorly by zero, for a diagonal matrix, and superiorly by 1.0, for a uniform matrix. Equation (B.1) provides a suitable *non-diagonality weighting function:* the higher the non-diagonality the higher the weight. Sparsification (noise-suppression) may be promoted by zeroing the weights above a cut-off frequency. According to our experience, such a weighting function generally allows satisfactory source estimation with continuously recorded EEG. We have observed that the non-diagonality function (B.1) is highly correlated with overall energy (trace of the cospectral matrices), but is not as much influenced by the dominant occipital rhythms (8-13Hz). In this fashion using a non-diagonality weighting function is in line with previous works in time-frequency BSS where the diagonalization effort has been concentrated on high-energy time-frequency regions (Belouchrani and Amin, 1998).

*C. Removing the DC-level (assuming zero-mean processes)*

Typically, BSS models assume zero-mean processes. Thus, for DC EEG amplifiers, the DC-level needs to be removed. Simply, the first cospectrum (0Hz) is not considered in the diagonalization set. Notice that FFT estimates at positive frequencies are not affected by the DC level (Bloomfield, 2000, p. 90), hence there is no need to remove the mean, detrend or band-pass the signal before computing the FFT (the same is not true for lagged covariance matrices).

*D. Evaluating the explained variance of source components*

In the introduction we have suggested that the energy of the output source can be evaluated in spite of their sign arbitrariness. For simplicity, we illustrate the method for diagonalization sets of the kind $\mathbb{C}:\{C_{(f)}\}$, that is, when only source coloration is exploited. The method readily extends to any number of indices. First, let us scale the rows of estimated separating matrix $\hat{B}$ so that they all have unit L2 norm. Because of the energy arbitrariness this operation does not alter the output of the BSS. Let $b_m^T$ and $a_m$ be, respectively, the normalized $m^{th}$ row of $\hat{B}$ (separating vector) and the $m^{th}$ column of $\hat{A} = \hat{B}^+$ (spatial pattern) associated with the $m^{th}$ component. Let $V \in \mathbb{R}^{N \cdot N}$ be the covariance matrix of the raw EEG data. Its diagonal elements $V_{nn}$ hold the variance (energy) of the $n^{th}$ EEG channel. Using (1.1) and (1.0) and ignoring the noise term in the latter, the *total explained variance* of the source components is given by

$$VAR_{TOT} = tr(\hat{A}\hat{B}V\hat{B}^T\hat{A}^T) \leq tr(V),$$

with strict equality if $M=N$, since then $\hat{A}\hat{B} = I$. Similarly, the *explained variance* of the $m^{th}$ component alone is such as

$$VAR_m = tr(a_m b_m^T V b_m a_m^T)$$

and we have $\sum_m VAR_m = VAR_{TOT}$. Notice that the explained variance or its relative portion $VAR_m/VAR_{TOT}$ can be evaluated for any discrete frequency using cospectral matrix $C_{(f)}$ instead of $V$ in computing both $VAR_{TOT}$ and $VAR_m$ above. In the same way, one may evaluate the explained variance for any frequency band pass region using instead the sum of cospectra within the region. This turns useful when we need to evaluate the energy of several components describing brain oscillations in a specific frequency band; for example, the several rhythms usually observed in the Alpha (8-12 Hz) range.

*E. Subspace reduction and pre-whitening*

When using many EEG sensors it may be useful to estimate fewer source components than sensors ($M<N$). Reducing the dimension of the input matrices makes them better conditioned, which enhance the performance of the separation (Pham and Cardoso, 2001)[6]. The *subspace reduction* may follow different strategies. For instance, one may use model-driven beamforming to attenuate the signal originating outside the region of interest (Rodríguez-Rivera et al., 2006; Congedo, 2006). Here we show how to perform subspace reduction to estimate the *M* most energetic source components, which is an extension of the common pre-whitening step. Let $C_{TOT} = \sum_v C_{(v)}$ be the sum of cospectral matrices forming the original (unreduced) diagonalization set. As in Eq. (1.7), *v* is a holder for a number of indeces. Now find a matrix $\boldsymbol{H} = [\boldsymbol{F} \; \boldsymbol{G}]^T \in \mathbb{R}^{N \cdot N}$, with the square brackets indicating matrix partition, such that $\boldsymbol{H}\boldsymbol{C}_{TOT}\boldsymbol{H}^T = \boldsymbol{I}$; $\boldsymbol{F} \in \mathbb{R}^{M \cdot N}$ holds the first $M$ rows of $\boldsymbol{H}$ (signal subspace) and $\boldsymbol{G} \in \mathbb{R}^{N-M \cdot N}$ the remaining rows (noise subspace). Note that for the case $\upsilon \triangleq f$, $\boldsymbol{C}_{TOT}$ is the sum of cospectral matrices at several frequencies and, if all Fourier frequencies are included in the diagonalization set, then matrix $\boldsymbol{H}$ is the well-known whitening matrix. For $\upsilon \triangleq ijk...$ our definition of $\boldsymbol{C}_{TOT}$ is the natural extension to obtain a "global" whitening matrix. Let us now factorize the separating matrix such as $\hat{\boldsymbol{B}} = \hat{\boldsymbol{E}}\boldsymbol{F}$, with $\hat{\boldsymbol{E}} \in \mathbb{R}^{M \cdot M}$. We obtain a new diagonalization set $D$ by applying the reduction to all cospectral matrices such as

$$D : \{\boldsymbol{D}_{(\upsilon)}\}, \boldsymbol{D}_{(\upsilon)} = \boldsymbol{F}\boldsymbol{C}_{(\upsilon)}\boldsymbol{F}^T.$$ The AJD problem (1.6) is now

$$\hat{\boldsymbol{E}} = AJD(D)$$

and we obtain the solution to the BSS problem as

$$\hat{\boldsymbol{B}} = \hat{\boldsymbol{E}}\boldsymbol{F} \in \mathbb{R}^{M \cdot N}.$$

---

[6] When the analyzed time interval is short the cospectral matrices may be non positive definite, a requirement of Pham's AJD algorithm. In this case the subspace reduction is necessary to obtain convergence. The algorithm of Ziehe et al. (2004) does not impose this restriction but does not allow explicit weighting.

Note that it is not necessary to constrain the AJD matrix $\boldsymbol{E}$ to be orthogonal, effectively circumventing the aforementioned drawback of pre-whitening the data. Finally, note that working with AJDC we do not need to compute the covariance matrix of the data at all.